\documentclass{amsart}
\usepackage{amssymb, latexsym}

\def\bi{\begin{itemize}}
\def\bs{\begin{split}}
\def\es{\end{split}}
\def\ba{\begin{align}}
\def\bas{\begin{align*}}
\def\ea{\end{align}}
\def\eas{\end{align*}}
\def\Im{{\operatorname{Im}}}

\def\Re{{\operatorname{Re}}}

\def\C{{\mathbb C}} % Jim changed to mathbb
\def\R{{{\mathbb R}}}
\def\Z{{{\mathbb Z}}}
\def\N{{{\mathbb N}}}

\def\emph#1{{\it #1}}
\def\textbf#1{{\bf #1}}

\newcommand{\ess}{{2+\varepsilon, \frac{2n(2+\varepsilon)}{n(2+\varepsilon)-4}}}

\newcommand{\hf}{L_t^{\infty}\dot H_x^{\frac12}}
\newcommand{\ep}{\frac{2p(2+\varepsilon)}{\varepsilon},\frac{np(2+\varepsilon)}{4+\varepsilon}}
\newcommand{\epp}{\frac{2p(2+\varepsilon)}{\varepsilon}, \frac{2np(2+\varepsilon)}{np(2+\varepsilon)-2\varepsilon}}
\newcommand{\eddp}{2,\frac{2n}{n+2}}

\newcommand{\er}{\frac{2(2+\varepsilon)}{\varepsilon}, \frac{n(2+\varepsilon)}{4+\varepsilon}}
\newcommand{\lnr}{\langle\nabla\rangle}
\newcommand{\lqlr}{{L^q_t L^r_x}}
\newcommand{\ir}{{I \times \R^n}}

\newcommand{\rr}{{\R \times \R^n}}

\newcommand{\llr}{{\lqlr (\rr)}}

\newcommand{\hk}{{\dot H^k}}

\newcommand{\nmt}{{\frac 4{n-2}}}

\newcommand{\eps}{{\varepsilon}}

\theoremstyle{plain}
\newtheorem{theorem}{Theorem}
\newtheorem{definition}[theorem]{Definition}
\newtheorem{proposition}[theorem]{Proposition}
\newtheorem{lemma}[theorem]{Lemma}

\theoremstyle{remark}
\newtheorem*{remark}{Remark}
\newtheorem*{remarks}{Remarks}

\numberwithin{equation}{section} \numberwithin{theorem}{section}

\begin{document}

\title[Global well-posedness and scattering for NLS]
{Global well-posedness and scattering for a class of nonlinear Schr\"odinger equations below the energy space}
\author{Monica Visan}
\address{University of California, Los Angeles}
\email{mvisan@math.ucla.edu}
\author{Xiaoyi Zhang}
\address{Academy of Mathematics and System Sciences, Chinese Academy of Sciences}
\email{xyzhang@amss.ac.cn}
\subjclass[2000]{35Q55} \keywords{Nonlinear Schr\"odinger equation, well-posedness}

\begin{abstract}
We prove global well-posedness and scattering for the nonlinear Schr\"odinger equation with power-type nonlinearity
\begin{equation*}
\begin{cases}%\label{eq}
i u_t +\Delta u = |u|^p u, \quad \frac{4}{n}<p<\frac{4}{n-2},\\
u(0,x) = u_0(x)\in H^s(\R^n), \quad n\geq 3,
\end{cases}
\end{equation*}
below the energy space, i.e., for $s<1$.  In \cite{ckstt:low7}, J.~Colliander, M.~Keel, G.~Staffilani, H.~Takaoka, and T.~Tao
established polynomial growth of the $H^s_x$-norm of the solution, and hence global well-posedness for initial data in $H^s_x$,
provided $s$ is sufficiently close to $1$. However, their bounds are insufficient to yield scattering.
In this paper, we use the \emph{a priori} interaction Morawetz inequality to show that scattering holds in $H^s(\R^n)$ whenever $s$ is larger
than some value $0<s_0(n,p)<1$.
\end{abstract}

\maketitle

\section{Introduction}

We study the initial value problem for the defocusing nonlinear Schr\"odinger equation
\begin{equation}\label{equation}
\begin{cases}
i u_t +\Delta u = |u|^p u\\
u(0,x) = u_0(x)\in H^s(\R^n),
\end{cases}
\end{equation}
where $u(t,x)$ is a complex-valued function in spacetime $\R\times\R^n$, $n\geq 3$, and the regularity $s$ is assumed to satisfy $0<s<1$.
Here, the value $p$ is assumed to be $L_x^2$-supercritical, i.e., $p>\tfrac{4}{n}$, but energy-subcritical, i.e., $p<\tfrac{4}{n-2}$.

This equation has Hamiltonian
\begin{equation}\label{energy}
E(u(t)):=\int_{\R^n} \Bigl[\tfrac{1}{2}|\nabla u (t,x)|^2 +\tfrac{1}{p+2}|u(t,x)|^{p+2}\Bigr]dx.
\end{equation}
As \eqref{energy} is preserved\footnote{To justify the energy conservation rigorously, one can approximate the data $u_0$ by smooth data,
and also approximate the nonlinearity by a smooth nonlinearity, to obtain a smooth approximate solution, obtain an energy conservation law
for that solution, and then take limits, using the local well-posedness theory.  We omit the standard details.  Similarly for the mass conservation
law and Morawetz type inequalities.} by the flow corresponding to \eqref{equation}, we shall refer to it as the \emph{energy} and often write $E(u)$
for $E(u(t))$.

A second conserved quantity we will rely on is the mass $\|u(t)\|^2_{L^2_x(\R^n)}$.

This equation has a natural scaling.  More precisely, the map
\begin{equation}\label{scaling}
u(t,x)\mapsto u_\lambda(t,x):=\lambda^{-\frac 2p}u\Bigl(\frac {t}{\lambda^2},\frac x{\lambda}\Bigr)
\end{equation}
maps a solution to \eqref{equation} to another solution to \eqref{equation}.
We define the critical regularity $s_c:=\frac n2-\frac 2p$; it is easy to verify that the scaling \eqref{scaling} leaves the $\dot H_x^{s_c}$-norm
invariant, up to a scaling of time.  In the case when $p=\frac{4}{n}$, $s_c=0$, which is why the nonlinearity $|u|^{\frac{4}{n}}u$ is called
$L_x^2$-critical.  When $p=\frac{4}{n-2}$, $s_c=1$, and hence the nonlinearity $|u|^{\nmt}u$ is called $\dot H^1_x$- or energy-critical.

The local and global theory for \eqref{equation} has been extensively studied.  It is known (see \cite{cwI}) that
the Cauchy problem \eqref{equation} is locally wellposed\footnote{By local well-posedness we mean existence, uniqueness, and uniform continuity
of the solution upon the initial data.} in $H^s_x$ for $s\ge \max\{0,s_c\}$.
%A global theory for small data is also available for the same regularity.
These results are known to be sharp in the sense that uniform continuity of the solution upon the initial data may fail
in the supercritical case $s<s_c$ (see \cite{cct:ill}).

By the local theory available in the subcritical case $s>s_c$ (specifically, the fact that the lifespan of the local solution depends only on the
$H_x^s$-norm of the initial data), global well-posedness for large data would follow immediately from a pointwise in time
$H^s_x$ bound on the solution and the usual iterative argument.  There are two cases when global well-posedness  follows from
the conservation laws for \eqref{equation}:  In the $L^2_x$-subcritical case ($s_c<0$), the conservation of mass implies global well-posedness in $L_x^2$.
In the energy-subcritical case ($s_c<1$), the Hamiltonian conservation combined with the Gagliardo-Nirenberg inequality yield
global well-posedness in $H^1_x$, both for the defocusing equation and for the focusing equation with $p<\frac 4n$.

This leaves open the question of global well-posedness in $H^s_x$ in the intermediate regime $0\le s_c \le s<1$.
The first result to address this problem belongs to Bourgain, \cite{Blowscatter}, who proved that the cubic defocusing NLS is globally wellposed
in $H^s(\R^3)$ for $s>\frac {11}{13}$.  Subsequently, J.~Colliander, M.~Keel, G.~Staffilani, H.~Takaoka, and T.~Tao
developed the `I-method', which they used to treat many problems including the one dimensional quintic NLS, and the two and three dimensional cubic NLS,
\cite{ckstt:low5, ckstt:low1, ckstt:low6, ckstt:low4, ckstt:low7, ckstt:low}.  For results in one dimension see also \cite{tzirakis}.

The idea behind the `I-method' is to smooth out the initial data (assumed to lie in $H_x^s$, $0<s<1$) in order to be able
to access the good local and global theory available at $H^1_x$ regularity.  To this end, one introduces the Fourier multiplier $I$,
which is the identity on low frequencies and behaves like a fractional integral operator of order $1-s$ on high frequencies.  Thus,
the operator $I$ maps $H_x^s$ to $H_x^1$.  However, even though we do have energy conservation for \eqref{equation}, $Iu$ is not a solution
to \eqref{equation} and hence, we expect an energy increment.  The key is to prove that $E(Iu)$ is an `almost conserved' quantity.
This requires delicate estimates on the commutator between $I$ and the nonlinearity.  When $p$ is an even integer, one can write the commutator
explicitly using the Fourier transform and control it by multilinear analysis and bilinear estimates; this type of estimates depend, of course,
on the exact form of the nonlinearity.  However, when $p$ is not an even integer, this method fails.

Relying on more rudimentary tools such as Taylor's expansion and Strichartz estimates, J.~Colliander, M.~Keel, G.~Staffilani, H.~Takaoka, and T.~Tao,
\cite{ckstt:low7}, proved global well-posedness for \eqref{equation} in $H^s(\R^n)$ for arbitrary $0<p<\frac{4}{n-2}$ and $s$ sufficiently close to $1$.
However, the polynomial upper bounds on the $H^s_x$-norm of the solution, which they obtain, are insufficient to derive scattering.

The goal of this paper is to prove that scattering does hold in $H^s_x$ for \eqref{equation}, for a sufficiently large regularity $s_0(n,p)<s<1$.
Here, the value $s_0(n,p)$ is defined by the following ghastly expression:
\begin{align}\label{ugly}
s_0(n,p)&=\max\{s_1(n,p), s_2(n,p), s_3(n,p)\},
\end{align}
where
\begin{align*}
s_1(n,p)&:=\tfrac{np}{2(p+2)}\\
s_2(n,p)&:=\tfrac{1+\min\{1,p\}s_c}{1+\min\{1,p\}}\\
s_3(n,p)&:=\min_{0<\sigma\leq \sigma_0} s_+(n,p,\sigma).
\end{align*}
Here, $\sigma_0$ must satisfy $2\sigma_0[8-p(n+2)]<(n-3)(pn-4)$ and $\sigma_0\leq s$, and $s_+(n,p,\sigma)$ is the larger of the two roots to
the quadratic equation
$$
s_c(1-s)[n-3-\sigma(n-6)]=\min\{1,p\}\sigma(s_c-s)^2.
$$
As $\sigma\mapsto s_+(n,p,\sigma)$ is a decreasing function (a straightforward but messy computation), $s_3(n,p)=s_+(n,p,\sigma_0)$.

\begin{remarks} \ \\
1. In dimensions four and higher, $s_2(n,p)$ can be omitted from  \eqref{ugly} because it is dominated by $s_1(n,p)$.  In dimension three,
$s_1(n,p)\leq s_2(n,p)$ only for $2\leq p<4$.

\noindent
2. In the case $n=3$, $s_+(n,p,\sigma)$ is independent of $\sigma$:
$$
s_3(n,p)=\frac{-s_c+ \sqrt{12s_c-3s_c^2}}{2}.
$$
From this and a little work, one discovers
$$
s_0(3,p)=
\begin{cases}
s_1(3,p), & \quad \text{for} \quad \frac{4}{3}<p<\frac{1+\sqrt{13}}{3}\\
s_3(3,p), & \quad \text{for} \quad \frac{1+\sqrt{13}}{3}<p<4.
\end{cases}
$$

\noindent
3. As $\sigma\searrow 0$, so $s_+(n,p,\sigma)\nearrow 1$.  More precisely, as $\sigma\searrow 0$, we have
$$
s_+(n,p,\sigma)=1-\frac{(s_c-1)^2\min\{1,p\}}{(n-3)s_c}\sigma + O(\sigma^2).
$$
Hence, the smaller $\sigma_0$ is (which we must require in order to treat values of $p$ close to the $L^2_x$-critical exponent,
see Lemma~\ref{interpolation1}) the larger we need to choose $s$.
\end{remarks}

The scattering theory in the energy class has been extensively studied.  J.~Ginibre and G.~Velo, \cite{gv:scatter}, proved scattering in $H^1_x$
for \eqref{equation} with $p$ in the long range, i.e., $\tfrac{4}{n}<p<\tfrac{4}{n-2}$.  Very recently, T.~Tao, M.~Visan, and X.~Zhang, \cite{tvz},
gave a new simpler proof of this result relying on the \emph{a priori} interaction Morawetz inequality.  In this paper, we will use the
interaction Morawetz estimate to prove scattering for \eqref{equation} below the energy space.

The \emph{a priori} interaction Morawetz inequality (for a proof in higher dimensions, see, for example, \cite{tvz}) gives
\begin{align}\label{a priori}
\int_I\int_{\R^{n}}\int_{\R^{n}} \frac{|u(t,y)|^{2} |u(t,x)|^{2}}{|x-y|^3} & dx\,dy\,dt
+\int_I\int_{\R^{n}}\int_{\R^{n}}\frac{|u(t,y)|^{2}|u(t,x)|^{p+2}}{|x-y|}dx\,dy\,dt \notag\\
& \lesssim \|u_0\|_{L_x^2}^2 \|u\|_{L_t^\infty \dot H^{\frac 12}_x(\ir)}^2
\end{align}
on any spacetime slab $\ir$ on which the solution $u$ to \eqref{equation} exists and lies in $H^{\frac 12}_x$.  A consequence of \eqref{a priori}
and some harmonic analysis (for a proof see again \cite{tvz}), is the following \emph{a priori} estimate on the solution to \eqref{equation}:
\begin{align}\label{negative deriv}
\||\nabla|^{-\frac{n-3}{4}} u\|_{L_{t,x}^4(\ir)}\lesssim \|u_0\|_{L_x^2}^{\frac 12} \|u\|_{L_t^\infty \dot H^{\frac 12}_x(\ir)}^{\frac 12}.
\end{align}
Interpolating between \eqref{negative deriv} and $u \in L_t^\infty \dot H^\sigma_x$ for $0<\sigma \leq s$, we obtain
\begin{equation*}
\|u\|_{M_\sigma(I\times\R^n)}
 \lesssim \bigl(\|u\|_{L_x^2}\|u\|_{\hf(\ir)}\bigr)^{\frac {2\sigma}{n-3+4\sigma}}\|u\|_{L_t^\infty\dot H^\sigma_x (\ir)}^{\frac{n-3}{n-3+4\sigma}},
\end{equation*}
where we define the Morawetz norm by
$$
\|u\|_{M_\sigma(\ir)}:=\|u\|_{L_t^{\frac{n-3+4\sigma}{\sigma}}L_x^{\frac{2(n-3+4\sigma)}{n-3+2\sigma}}(\ir)}.
$$

Our main result is the following
\begin{theorem}\label{scattering}
Let $n\ge 3$, $s>s_0(n,p)$ and let $u_0\in H^s(\R^n)$.  Then, the Cauchy problem \eqref{equation} is globally wellposed and the global solution $u$
enjoys the following uniform bound
\begin{align*}
\|u\|_{L_t^\infty H^s_x(\R\times \R^n)} \le C(\|u_0\|_{H^s_x}).
\end{align*}
Furthermore, there exist unique scattering states $u_\pm\in H^s_x$ such that
$$
\|u(t)-e^{it\Delta}u_\pm\|_{H^s_x}\to 0 \quad \text{as} \quad t\to \pm \infty.
$$
\end{theorem}

We record part of our results and compare them with the best known results in the table below.  As the reader can see, our results are not optimal
when $p$ is an even integer.  However, our method here is robust and does not depend on the exact form of nonlinearity.

\begin{table}[h]
\begin{center}
\def\up{\vrule width 0mm height 2.5ex depth 1.2ex}
\begin{tabular}{|c|c|c|c|l|}
\hline $n$   & $p$          & $s_c$          & $s$  (best known)              & $s$ (ours)\\
\hline \up$3$   & $2$          & $\tfrac 12$    & $0.8$ (\cite{ckstt:low})       & 0.895 \\
\hline \up$3$   & $3$          & $\tfrac 56$    & --                             & 0.990\\
\hline \up$4$   & $\tfrac 32$  & $\tfrac 23$    & --                             & 0.958\\
\hline
\end{tabular}
\end{center}
\label{Table1}
\end{table}

The remaining part of the paper is organized as follows:  In Section~2 we introduce notation and prove some lemmas that will be useful.
In Section~3 we prove Theorem~\ref{scattering}.

%%%%%%%%%%%%%%%%%%%%%%%%%%%%%%%%%%%%%%%%%%%%%%%%%%%%%%%%%%%%%%%%%%%%%%%%%%%%%%%%%%%%%%%%%%%
%
%
%                                   Section
%
%
%%%%%%%%%%%%%%%%%%%%%%%%%%%%%%%%%%%%%%%%%%%%%%%%%%%%%%%%%%%%%%%%%%%%%%%%%%%%%%%%%%%%%%%%%%%

\section{Preliminaries}
We will often use the notation $X \lesssim Y$ whenever there exists some constant $C$ so that $X \leq CY$. Similarly, we will use $X
\sim Y$ if $X \lesssim Y \lesssim X$.  We use $X \ll Y$ if $X \leq cY$ for some small constant $c$. The derivative operator $\nabla$
refers to the space variable only.  We use $A\pm$ to denote $A\pm\eps$ for any sufficiently small $\eps>0$.

Let $F(z):=|z|^pz$ be the function that defines the nonlinearity in \eqref{equation}.  Then,
$$
F_z(z):=\frac{\partial F}{\partial z}(z)=\frac{p+2}2|z|^{p} \quad \text{and}
\quad F_{\bar z}(z):=\frac{\partial F}{\partial \bar z}(z)=\frac{p}2|z|^{p}\frac{z}{\bar z}.
$$
We write $F'$ for the vector $(F_z,F_{\bar z})$ and adopt the notation
$$
w\cdot F'(z):=wF_{z}(z)+\bar wF_{\bar z}(z).
$$
In particular, we observe the chain rule
$$
\nabla F(u)=\nabla u\cdot F'(u).
$$
Clearly $F'(z)=O(|z|^{p})$ and we have the H\"older continuity estimate
\begin{align}\label{holder continuity}
|F'(z)-F'(w)|\lesssim |z-w|^{\min\{1,p\}}(|z|+|w|)^{p-\min\{1,p\}}
\end{align}
for all $z, w\in \C$.  By the Fundamental Theorem of Calculus,
$$
F(z+w)-F(z)=\int_0^1 w\cdot F'(z+\theta w)d\theta
$$
and hence
$$
F(z+w)=F(z)+O(|w||z|^p)+O(|w|^{p+1})
$$
for all complex values $z$ and $w$.

We use $L_x^r(\R^n)$ to denote the Banach space of functions $f:\R^n\to \C$ whose norm
$$
\|f\|_r:=\Bigl(\int_{\R^n} |f(x)|^r dx\Bigr)^{\frac{1}{r}}
$$
is finite, with the usual modifications when $r=\infty$.  For any
non-negative integer $k$, we denote by $H^{k,r}(\R^n)$ the Sobolev
space defined as the closure of test functions in the norm
$$
\|f\|_{H^{k,r}}:=\sum_{|\alpha|\leq k}\Bigl\|\frac{\partial
^\alpha}{\partial x^\alpha}f\Bigr\|_r.
$$
We will often denote $H^{k,2}$ by $H^k$.

We use $\lqlr$ to denote the spacetime norm
$$
\|u\|_{q,r}:=\|u\|_{\llr} :=\Bigl(\int_{\R}\Bigl(\int_{\R^n} |u(t,x)|^r dx \Bigr)^{q/r} dt\Bigr)^{1/q},
$$
with the usual modifications when either $q$ or $r$ are infinity, or when the domain $\R \times \R^n$ is
replaced by some smaller spacetime region.  When $q=r$ we abbreviate $\lqlr$ by $L^q_{t,x}$.

We define the Fourier transform on $\R^n$ to be
$$
\hat f(\xi) := \int_{\R^n} e^{-2 \pi i x \cdot \xi} f(x) dx.
$$

We will make use of the fractional differentiation operators $|\nabla|^s$ defined by
$$
\widehat{|\nabla|^sf}(\xi) := |\xi|^s \hat f (\xi).
$$
These define the homogeneous Sobolev norms
$$
\|f\|_{\dot H^{s}_x} := \| |\nabla|^s f \|_{L^2_x}.
$$

Let $e^{it\Delta}$ be the free Schr\"odinger propagator.  In physical space this is given by the formula
$$
e^{it\Delta}f(x) = \frac{1}{(4 \pi i t)^{n/2}} \int_{\R^n} e^{i|x-y|^2/4t} f(y) dy
$$
for $t\neq 0$ (using a suitable branch cut to define $(4\pi it)^{n/2}$),
while in frequency space one can write this as
\begin{equation}\label{fourier rep}
\widehat{e^{it\Delta}f}(\xi) = e^{-4 \pi^2 i t |\xi|^2}\hat f(\xi).
\end{equation}
In particular, the propagator obeys the \emph{dispersive inequality}
\begin{equation}\label{dispersive ineq}
\|e^{it\Delta}f\|_{L^\infty_x} \lesssim
|t|^{-\frac{n}{2}}\|f\|_{L^1_x}
\end{equation}
for all times $t\neq 0$.

We also recall \emph{Duhamel's formula}
\begin{align}\label{duhamel}
u(t) = e^{i(t-t_0)\Delta}u(t_0) - i \int_{t_0}^t e^{i(t-s)\Delta}(iu_t + \Delta u)(s) ds.
\end{align}

\begin{definition}
A pair of exponents $(q,r)$ is called Schr\"odinger-\emph{admissible} if
$$
\frac{2}{q} +\frac{n}{r} = \frac{n}{2},\quad 2 \leq q,r \leq \infty, \quad \text{and} \quad (q,r,n)\neq (2,\infty,2).
$$
\end{definition}
Throughout this paper we will use the following admissible pairs:
$$
(2,\tfrac{2n}{n-2}) \quad \text{and} \quad (\tfrac{n-3+4\sigma}{\sigma}, \tfrac{2n(n-3+4\sigma)}{n(n-3+4\sigma)-4\sigma}) \quad \text{for} \, \sigma>0
$$
$$
(\tfrac{2p(2+\eps)}{\eps}, \tfrac{2np(2+\eps)}{np(2+\eps)-2\eps}) \quad \text{and} \quad (2+\eps, \tfrac{2n(2+\eps)}{n(2+\eps)-4}) \quad \text{for some small} \ \eps>0.
$$

We record the standard Strichartz estimates which we will invoke repeatedly throughout this paper (for a proof see \cite{tao:keel}):

\begin{lemma}\label{lemma linear strichartz}
Let $I$ be a compact time interval, $t_0\in I$, $k$ an arbitrary integer, and let $u$ be a solution to the forced Schr\"odinger equation
\begin{equation*}
i u_t + \Delta u =\sum_{i=1}^m F_i
\end{equation*}
for some functions $F_1$, $\cdots$, $F_m$.  Then,
\begin{equation}
\||\nabla| ^k u\|_{L^q_tL_x^r(I\times\R^n)} \lesssim \|u(t_0)\|_{\hk(\R^n)} +  \sum_{i=1}^m \||\nabla| ^k F_i\|_{L_t^{q_i'}L_x^{r_i'} (I\times\R^n)}
\end{equation}
for any admissible pairs $(q,r)$ and $(q_i,r_i)$, $1\leq i\leq m$.
\end{lemma}

We will also need some Littlewood-Paley theory.  Specifically, let $\varphi(\xi)$ be a smooth bump supported in the ball
$|\xi| \leq 2$ and equalling one on the ball $|\xi| \leq 1$.  For each dyadic number $N \in 2^\Z$ we define the
Littlewood-Paley operators
\begin{align*}
\widehat{P_{\leq N}f}(\xi) &:=  \varphi(\xi/N)\hat f (\xi),\\
\widehat{P_{> N}f}(\xi) &:=  [1-\varphi(\xi/N)]\hat f (\xi),\\
\widehat{P_N f}(\xi) &:=  [\varphi(\xi/N) - \varphi (2 \xi /N)] \hat f (\xi).
\end{align*}
Similarly we can define $P_{<N}$, $P_{\geq N}$, and $P_{M < \cdot \leq N} := P_{\leq N} - P_{\leq M}$, whenever $M$ and
$N$ are dyadic numbers.  We will frequently write $f_{\leq N}$ for $P_{\leq N} f$ and similarly for the other operators.
We recall the following standard Bernstein and Sobolev type inequalities:

\begin{lemma}\label{bernstein}
For any $1\le p\le q\le\infty$ and $s>0$, we have
\begin{align*}
\|P_{\geq N} f\|_{L^p_x} &\lesssim N^{-s} \| |\nabla|^s P_{\geq N} f \|_{L^p_x}\\
\| |\nabla|^s  P_{\leq N} f\|_{L^p_x} &\lesssim N^{s} \| P_{\leq N} f\|_{L^p_x}\\
\| |\nabla|^{\pm s} P_N f\|_{L^p_x} &\sim N^{\pm s} \| P_N f \|_{L^p_x}\\
\|P_{\leq N} f\|_{L^q_x} &\lesssim N^{\frac{n}{p}-\frac{n}{q}} \|P_{\leq N} f\|_{L^p_x}\\
\|P_N f\|_{L^q_x} &\lesssim N^{\frac{n}{p}-\frac{n}{q}} \| P_N f\|_{L^p_x}.
\end{align*}
\end{lemma}

\vspace{0.4cm}

For $N>1$, we define the Fourier multiplier $I:=I_N$ by
$$
\widehat{I_N u}(\xi):=m_N(\xi)\hat u(\xi),
$$
where $m_N$ is a smooth radial decreasing cutoff function such that
$$
m_N(\xi)=\left\{
\begin{array}{cc}
1, \quad \text{if} \quad  |\xi|\le N\\
\bigl(\frac{|\xi|}{N}\bigr)^{s-1}, \quad \text{if} \quad |\xi|\ge 2N.
\end{array}
\right.
$$
Thus, $I$ is the identity operator on frequencies $|\xi|\le N$ and behaves like a fractional integral
operator of order 1-s on higher frequencies. In particular, $I$ maps $H^s_x$ to $H_x^1$; this allows us to
access the good local and global theory available for $H_x^1$ data.  We collect the basic properties of $I$ into the following:

\begin{lemma}\label{basic property}
Let $1<p<\infty$ and $0\leq\sigma\le s<1$.  Then,
\begin{align}
\|If\|_p&\lesssim \|f\|_p \label{i1}\\
\||\nabla|^\sigma P_{> N}f\|_p&\lesssim N^{\sigma-1}\|\nabla I f\|_p \label{i2}\\
\|f\|_{H^s_x}\lesssim \|If\|_{H^1_x}&\lesssim N^{1-s}\|f\|_{H^s_x}.\label{i3}
\end{align}
\end{lemma}
\begin{proof}
The estimate \eqref{i1} is a direct consequence of the multiplier theorem.

To prove \eqref{i2}, we write
$$
\||\nabla|^\sigma P_{> N} f\|_p=\|P_{> N}|\nabla |^\sigma(\nabla I)^{-1}\nabla I f\|_p.
$$
The claim follows again from the multiplier theorem.

Now we turn to \eqref{i3}.  By the definition of the operator $I$ and \eqref{i2},
\begin{align*}
\|f\|_{H^s_x}&\lesssim \|P_{\le N} f\|_{H^s_x}+\|P_{>N}f\|_2+\||\nabla|^s P_{>N} f\|_2\\
&\lesssim \|P_{\le N} I f\|_{H_x^1}+N^{-1}\|\nabla I f\|_2+N^{s-1}\|\nabla I f\|_2\\
&\lesssim \|If\|_{H^1_x}.
\end{align*}
On the other hand, since the operator $I$ commutes with $\langle \nabla\rangle^s$,
\begin{align*}
\|If\|_{H_x^1}
=\|\langle\nabla\rangle^{1-s}I\langle\nabla\rangle^s f\|_2
\lesssim N^{1-s}\|\langle\nabla\rangle^sf\|_2
\lesssim N^{1-s}\|f\|_{H^s_x},
\end{align*}
which proves the last inequality in \eqref{i3}.  Note that a similar argument also yields
\begin{align}\label{i4}
\|If\|_{\dot H^1_x}&\lesssim N^{1-s}\|f\|_{\dot H^s_x}.
\end{align}
\end{proof}

The estimate \eqref{i2} shows that we can control the high
frequencies of a function $f$ in the Sobolev space $H^{\sigma,p}$ by
the smoother function $If$ in a space with a loss of derivative but
a gain of negative power of $N$. This fact is crucial in extracting
the negative power of $N$ when estimating the increment of the
modified Hamiltonian.

When $p$ is an even integer, one can use multilinear analysis to understand commutator expressions like $F(Iu) - IF(u)$; on the Fourier side,
one can expand this commutator into a product of Fourier transforms of $u$ and $Iu$ and carefully measure the frequency interactions to derive
an estimate (see for example \cite{ckstt:low}).  However, this is not possible when $p$ is not an even integer.  Instead, we will have to rely
on the following rougher (weaker, but more robust) lemma:

\begin{lemma}\label{bilinear}
Let $1<r,r_1,r_2<\infty$ be such that $\frac 1 r=\frac 1{r_1}+\frac 1{r_2}$ and let $0<\nu<s$.  Then,
\begin{align}\label{bilinear eq}
\|I(fg)-(If)g\|_r\lesssim N^{-(1-s+\nu)}\|If\|_{r_1}\|\langle\nabla\rangle^{1-s+\nu}g\|_{r_2}.
\end{align}
\end{lemma}
\begin{proof} Applying a Littlewood-Paley decomposition to $f$ and $g$, we write
\begin{align}
I(fg)-(If)g
&=I(fg_{\le 1})-(If)g_{\le1} +\sum_{1<M\in 2^{\Z}}\bigl[I(f_{\lesssim M}g_M)-(If_{\lesssim M})g_M \bigr] \nonumber\\
&\qquad \qquad +\sum_{1<M\in 2^{\Z}}\bigl[I(f_{\gg M}g_M)-(If_{\gg M})g_M\bigr] \nonumber\\
&=I(f_{\gtrsim N}g_{\le 1})-(If_{\gtrsim N}) g_{\le 1} +\sum_{N\lesssim M\in 2^{\Z}}\bigl[I(f_{\lesssim M} g_M) -(If_{\lesssim M})g_M\bigr] \nonumber\\
&\qquad\qquad +\sum_{1<M\in 2^{\Z}}\bigl[I(f_{\gg M}g_M)-(If_{\gg M})g_M \bigr]\nonumber\\
&=I +II +III.\label{difference}
\end{align}
The second equality above follows from the fact that the operator $I$ is the identity operator on frequencies $|\xi|\leq N$; thus,
\begin{align*}
I(f_{\ll N}g_{\le1})=(If_{\ll N})g_{\le 1} \quad \text{and} \quad I(f_{\lesssim M}g_{M})=(If_{\lesssim M})g_{M}& \quad \text{for all} \quad M\ll N.
\end{align*}

We first consider $II$.  Dropping the operator $I$, by H\"older and Bernstein we estimate
\begin{align*}
\|I(f_{\lesssim M} g_M) -(If_{\lesssim M})g_M\|_r&\lesssim \|f_{\lesssim M}\|_{r_1}\|g_M\|_{r_2}\\
&\lesssim \bigl(\frac MN\bigr)^{1-s} \|If\|_{r_1}\|g_M\|_{r_2}\\
&\lesssim M^{-\nu}N^{-(1-s)}\|If\|_{r_1}\||\nabla|^{1-s+\nu}g\|_{r_2}.
\end{align*}
Summing over all $N\lesssim M\in 2^{\Z}$, we get
\begin{equation}\label{d1}
II\lesssim N^{-(1-s+\nu)}\|If\|_{r_1}\||\nabla|^{1-s+\nu}g\|_{r_2}.
\end{equation}

We turn now towards $III$.  Applying a Littlewood-Paley decomposition to $f$, we write each term in $III$ as
\begin{align*}
I(f_{\gg M}g_M)-(If_{\gg M})g_M
&=\sum_{1\ll k\in\N} \bigl[ I(f_{2^k M}g_M)-(If_{2^kM})g_M\bigr]\\
&=\sum_{\substack{1\ll k\in\N \\ N\lesssim 2^kM}} \bigl[I(f_{2^k M} g_M)-(If_{2^kM})g_M\bigr].
\end{align*}
To derive the second inequality, we used again the fact that the operator $I$ is the identity on frequencies $|\xi|\le N$.

We write
$$
[I(f_{2^k M} g_M)-(If_{2^kM})g_M]\widehat{\ }(\xi)=\int_{\xi=\xi_1+\xi_2}(m_N(\xi_1+\xi_2)-m_N(\xi_1))\widehat{f_{2^k M}}(\xi_1)\widehat{g_M}(\xi_2).
$$
For $|\xi_1|\sim 2^kM$, $k\gg 1$, and $|\xi_2|\sim M$, the Fundamental Theorem of Calculus implies
$$
|m_N(\xi_1+\xi_2)-m_N(\xi_1)|\lesssim 2^{-k} \bigl(\frac{2^kM}{N}\bigr)^{s-1}.
$$
By the Coifman-Meyer multilinear multiplier theorem, \cite{coifmey:1, coifmey:2}, and Bernstein, we get
\begin{align*}
\|I(f_{2^k M} g_M)-(If_{2^k M})g_M\|_r
&\lesssim 2^{-k} \bigl(\frac{2^kM}{N}\bigr)^{s-1}\|f_{2^k M}\|_{r_1}\|g_M\|_{r_2}\\
&\lesssim 2^{-k}M^{-(1-s+\nu)}\|If\|_{r_1}\||\nabla|^{1-s+\nu} g\|_{r_2}.
\end{align*}
Summing over $M$ and $k$ such that $N\lesssim 2^kM$, and recalling that $0<\nu<s$, we get
\begin{align}\label{d2}
III\lesssim N^{-(1-s+\nu)}\|If\|_{r_1}\||\nabla|^{1-s+\nu}g\|_{r_2}.
\end{align}

The estimate $I$, we apply the same argument as for $III$.  We get
\begin{align}\label{d3}
I=\|I(f_{\gtrsim N}g_{\le 1})-(If_{\gtrsim N})g_{\le 1} \|_r
&\lesssim \sum_{k\in \mathbb N,2^k \gtrsim N}\|I(f_{2^k}g_{\le 1})-(If_{2^k})g_{\le 1}\|_r \nonumber\\
&\lesssim \sum_{k\in \mathbb N,2^k\gtrsim N} 2^{-k}\|If\|_{r_1}\|g\|_{r_2}\nonumber\\
&\lesssim N^{-1}\|If\|_{r_1}\|g\|_{r_2}.
\end{align}
Putting \eqref{difference} through \eqref{d3} together, we derive \eqref{bilinear eq}.
\end{proof}

As an application of Lemma \ref{bilinear} we have the following commutator estimate:

\begin{lemma}\label{nablaif}
Let $1<r,r_1,r_2<\infty$ be such that $\frac 1r=\frac 1{r_1}+\frac 1{r_2}$.  Then, for any $0<\nu<s$ we have
\begin{align}
\|\nabla IF(u)-(I\nabla u)F'(u)\|_r\lesssim N^{-1+s-\nu}\|\nabla I u\|_{r_1}\|\langle\nabla\rangle^{1-s+\nu}F'(u)\|_{r_2}\label{na1}\\
\|\nabla I F(u)\|_r\lesssim \|\nabla I u\|_{r_1}\|F'(u)\|_{r_2}+N^{-1+s-\nu} \|\nabla Iu\|_{r_1}\|\langle\nabla\rangle^{1-s+\nu}F'(u)\|_{r_2}. \label{na2}
\end{align}
\end{lemma}
\begin{proof}
As
$$
\nabla F(u)=F'(u)\nabla u,
$$
the estimate \eqref{na1} follows immediately from Lemma \ref{bilinear} with $f:=\nabla u$ and $g:=F'(u)$. The estimate \eqref{na2} is a
consequence of \eqref{na1} and the triangle inequality.
\end{proof}

Since we work at regularity $0<s<1$, we will need the following fractional chain rule to estimate our nonlinearity in $H_x^s$.

\begin{lemma}[Fractional chain rule for a $C^1$ function, \cite{chris:weinstein}]\label{F Lip}
Suppose that $F\in C^1(\mathbb C)$, $\alpha \in (0,1)$, and $1<r,r_1,r_2<\infty$ such that $\frac 1r=\frac 1{r_1}+\frac 1{r_2}$.   Then,
$$
\||\nabla|^{\alpha}F(u)\|_r\lesssim \|F'(u)\|_{r_1}\||\nabla|^{\alpha}u\|_{r_2}.
$$
\end{lemma}

When the function $F$ is no longer $C^1$, but merely H\"older continuous, we have the following useful chain rule:

\begin{lemma}[\cite{Monica:thesis}]\label{fdfp}
Let $F$ be a H\"older continuous function of order $0<\alpha<1$.  Then, for every $0<\sigma<\alpha$, $1<r<\infty$,
and $\tfrac{\sigma}{\alpha}<\rho<1$ we have
\begin{align}\label{fdfp2}
\bigl\| |\nabla|^\sigma F(u)\bigr\|_r
\lesssim \bigl\||u|^{\alpha-\frac{\sigma}{\rho}}\bigr\|_{r_1} \bigl\||\nabla|^\rho u\bigr\|^{\frac{\sigma}{\rho}}_{\frac{\sigma}{\rho}r_2},
\end{align}
provided $\tfrac{1}{r}=\tfrac{1}{r_1} +\tfrac{1}{r_2}$ and $(1-\frac\sigma{\alpha \rho})r_1>1$.
\end{lemma}

In Section~3, we will need to control the nonlinearity in terms of the Morawetz norm $M_\sigma$.
The idea is simple.  Notice that by Lemma~\ref{F Lip} and H\"older's inequality, we have
$$
\|\langle\nabla\rangle^s(|u|^p u)\|_{2,\frac {2n}{n+2}}\lesssim \|\langle\nabla\rangle^su\|_{2,\frac {2n}{n-2}}\|u\|^p_{\infty,\frac {np}2}.
$$
In order to get a factor of $\|u\|_{M_\sigma}$ on the right-hand side, we replace the space $L_t^\infty L_x^{\frac{np}{2}}$
by a space which interpolates between $L_t^{\frac{n-3+4\sigma}{\sigma}} L_x^{\frac{2(n-3+4\sigma)}{n-3+2\sigma}}$, $L_t^\infty L_x^{\frac{2n}{n-2s}}$
(which imbeds in $\dot H^s_x$), and the mass, $L_t^\infty L_x^2$.
To this end, we replace the space $L_t^2L_x^{\frac{2n}{n-2}}$ by the Strichartz space $L_t^{2+\eps} L_x^{\frac {2n(2+\eps)}{n(2+\eps)-4}}$
(for a sufficiently small $\eps>0$).  More precisely, we have

\begin{lemma}\label{use Morawetz}
Let $0<\sigma\leq s<1$ such that $\tfrac{\sigma(n-2)}{n-3+4\sigma}<s$ and let $\tfrac{4}{n}<p<\tfrac{4}{n-2s}$.  Then, there exists $\eps>0$
sufficiently small such that on every slab $\ir$ we have
\begin{align}\label{use Morawetz eq}
\|\langle\nabla\rangle^s(|u|^pu)\|_{2,\frac {2n}{n+2}}
\lesssim \|\langle\nabla\rangle^s u\|_{2+\eps,\frac {2n(2+\eps)}{n(2+\eps)-4}}\|u\|^{\frac{\eps(n-3+4\sigma)}{2\sigma(2+\eps)}}_{M_\sigma} \|u\|_{L_t^\infty L_x^2}^{\alpha(\eps)}\|u\|^{\beta(\eps)}_{L_t^\infty \dot H^s_x}.
\end{align}
Here,
$$
\alpha(\eps):=p\bigl(1-\tfrac{n}{2s}\bigr)+\tfrac{8\sigma+\eps[\sigma(n+2)-s(n-3+4\sigma)]}{2s\sigma(2+\eps)}   \quad \text{and} \quad \beta(\eps):=\tfrac{n}{s}\bigl(\tfrac{p}{2}-\tfrac{8+\eps(n+2)}{2n(2+\eps)}\bigr).
$$
\end{lemma}

\begin{proof}
Once $\alpha(\eps)$ and $\beta(\eps)$ are positive, the estimate \eqref{use Morawetz eq} is a direct consequence of H\"older's inequality and Sobolev
embedding, as the reader can easily check.  It is not hard to see that $\eps\mapsto \alpha(\eps)$ and $\eps\mapsto \beta(\eps)$ are decreasing functions
(for the former we need $s>\tfrac{\sigma(n-2)}{n-3+4\sigma}$).  Moreover,
$$
\alpha(\eps)\to p\bigl(1-\tfrac{n}{2s}\bigr) +\tfrac{2}{s} \quad \text{and} \quad \beta(\eps)\to \tfrac{n}{s}\bigl(\tfrac{p}{2}-\tfrac{2}{n}\bigr)
\quad \text{as} \quad \eps\to 0.
$$
As $\tfrac{4}{n}<p<\tfrac{4}{n-2s}$, the two limits are positive.  Thus, for $\eps>0$ sufficiently small we obtain
$$
\alpha(\eps)>0 \quad \text{and} \quad \beta(\eps)>0.
$$
This concludes the proof of Lemma~\ref{use Morawetz}.
\end{proof}

\begin{remark}
Note that the function $(0,s]\ni\sigma\mapsto\tfrac{\sigma(n-2)}{n-3+4\sigma}$ is increasing and attains its largest value for $\sigma=s$.  In this case,
the condition $\tfrac{\sigma(n-2)}{n-3+4\sigma}<s$ becomes $s>\tfrac{1}{4}$, which is implied by $s>s_1(n,p)=\tfrac{np}{2(p+2)}$.
\end{remark}

Finally, we will need the following

\begin{lemma}\label{interpolation1}
Let $0<\sigma\leq s<1$ and $\tfrac{4(n-3+4\sigma)}{n(n-3+2\sigma)+4\sigma}<p<\tfrac{4}{n-2s}$.  Then
$$
\|u\|_{\ep}\lesssim\|u\|_{M_\sigma}^{\theta}\|\nabla u\|_{\tfrac{n-3+4\sigma}{\sigma}, \tfrac{2n(n-3+4\sigma)}{n(n-3+4\sigma)-4\sigma}}^{1-\theta},
$$
where
$$
\eps:=\tfrac{4p\sigma}{n-3+4\sigma-2p\sigma} \quad \text{and} \quad \theta:=\tfrac{(n-3+4\sigma)[4-p(n-2)]}{2p[n-3-\sigma(n-6)]}.
$$
\end{lemma}
The proof of Lemma~\ref{interpolation1} involves straightforward computations using H\"older's inequality, interpolation, and Sobolev embedding
and we omit it.  Note that
$$
\tfrac{4(n-3+4\sigma)}{n(n-3+2\sigma)+4\sigma}\to \tfrac{4}{n} \quad \text{as}\quad \sigma\to 0
$$
and hence we can treat any $p\in (\tfrac{4}{n},\tfrac{4}{n-2s})$.  However, as $p$ approaches the $L_x^2$-critical value, i.e., $\tfrac{4}{n}$,
we are forced to choose $\sigma$ very small, which in turn forces $s$ to be close to $1$.

%%%%%%%%%%%%%%%%%%%%%%%%%%%%%%%%%%%%%%%%%%%%%%%%%%%%%%%%%%%%%%%%%%%%%%%%%%%%%%%%%%%%%%%%%%%
%
%
%                                   Section
%
%
%%%%%%%%%%%%%%%%%%%%%%%%%%%%%%%%%%%%%%%%%%%%%%%%%%%%%%%%%%%%%%%%%%%%%%%%%%%%%%%%%%%%%%%%%%%

\section{Proof of Theorem~\ref{scattering}}
\subsection{Global well-posedness}

In the local well-posedness theory, the time of existence of the unique solution to \eqref{equation} depends only on the $H_x^s$-norm
of the initial data.  Thus, by the usual iterative argument, global well-posedness would follow from a global $L_t^\infty H_x^s$ bound on the solution.

However, the $H^s_x$-norm of the solution is not a conserved quantity.  Nevertheless, it can be controlled by the $H_x^1$-norm
of the modified solution $Iu$ (see \eqref{i3}).  While we do have conservation of energy for \eqref{equation}, $Iu$ is not a solution
to \eqref{equation} and hence we expect an energy increment.  This will be proved to be small on intervals where the Morawetz norm is small,
which transfers the problem to controlling the Morawetz norm globally.  This idea is encapsulated in the following statement, which is proved
at the end of the section.

\begin{proposition}\label{energy increment scattering}
Let $s_0(n,p)<s<1$ and $\tfrac{4}{n}<p<\tfrac 4{n-2s}$.  Let $u$ be an $H^s_x$ solution to \eqref{equation} on $[t_0,T]$ with
$\|\nabla I u(t_0)\|_2\leq 1$.  Suppose in addition that
$$
\|u\|_{M_\sigma([t_0,T]\times \R^n)}\leq\eta
$$
for a sufficiently small $\eta>0$.  Then, for $N$ sufficiently large,
$$
\sup_{t\in[t_0,T]}E(Iu(t))=E(Iu(t_0))+O(N^{\min\{1,p\}(s_c-s)+}).
$$
Here, the implicit constant depends only on the size of $E(Iu(t_0))$.
\end{proposition}

Therefore, the proof of global well-posedness has been reduced to showing
\begin{align}\label{global Morawetz}
\|u\|_{M_\sigma(\R\times\R^n)}\leq C(\|u_0\|_{H_x^s}).
\end{align}
This also implies scattering, as we will show below.

Recall that interpolating between the \emph{a priori} interaction Morawetz inequality \eqref{negative deriv} and $L_t^\infty \dot H_x^\sigma$,
$0<\sigma\leq s$, we get
\begin{equation}\label{inter Mora esti}
\|u\|_{M\sigma(I\times\R^n)}\lesssim \bigl(\|u\|_{L_x^2}\|u\|_{\hf(\ir)}\bigr)^{\frac {2\sigma}{n-3+4\sigma}}\|u\|_{L_t^\infty\dot H^\sigma_x (\ir)}^{\frac{n-3}{n-3+4\sigma}}
\end{equation}
on any spacetime slab $\ir$ on which the solution to \eqref{equation} exists and lies in $H_x^{\max\{\frac 12,\sigma\}}$.
However, the $H_x^{\max\{\frac 12,\sigma\}}$-norm of the solution is not conserved either and in order to control it, we must resort
to the $H^s_x$ bound on the solution.  Thus, in order to obtain a global Morawetz estimate we need a global $H^s_x$ bound.
This sets us up for a bootstrap argument.

Let $u(x,t)$ be the solution to \eqref{equation}.  As $E(Iu_0)$ is not necessarily small, we will rescale the solution such that the energy of the
rescaled initial data satisfies the hypothesis of Proposition~\ref{energy increment scattering}.  Indeed, by scaling,
$$
u^{\lambda}(x,t):=\lambda^{-\frac 2p}u\Bigl(\frac x{\lambda},\frac t{\lambda^2}\Bigr)
$$
is a solution to \eqref{equation} with initial data
$$
u_0^{\lambda}:=\lambda^{-\frac 2p}u_0\Bigl(\frac x{\lambda}\Bigr).
$$
By \eqref{i4} and Sobolev embedding (we need $s>\frac{np}{2(p+2)}$),
\begin{align*}
&\|\nabla Iu_0^\lambda\|_2\lesssim N^{1-s}\|u_0^\lambda\|_{\dot H^s_x}=N^{1-s}\lambda^{s_c-s} \|u_0\|_{\dot H^s_x}\\
&\|u_0^\lambda\|_{p+2}=\lambda^{-\frac{2}{p}+\frac{n}{p+2}}\|u_0\|_{p+2}\lesssim\lambda^{-\frac{2}{p}+\frac{n}{p+2}}\|u_0\|_{H_x^s}.
\end{align*}
As we are in the energy subcritical case, $-\frac{2}{p}+\frac{n}{p+2}<0$.  Thus, taking $\lambda$ sufficiently large depending on $\|u_0\|_{H_x^s}$
and $N$ (which will be chosen later and will depend only on $\|u_0\|_{H_x^s}$), we get
$$
E(Iu_0^\lambda)\leq c(\|u_0\|_{H^s_x})  \ll 1.
$$

We now show that there exists an absolute constant $C_1$ such that
\begin{align}\label{rescaled Mora}
\|u^{\lambda}\|_{M_\sigma(\R\times\R^n)}\le C_1\lambda^{\frac {s_c[n-3-\sigma(n-6)]}{n-3+4\sigma}}.
\end{align}
Undoing the scaling, this yields \eqref{global Morawetz}.

By time reversal symmetry, it suffices to argue for positive times only.  Define
\begin{equation*}
\Omega_1:=\{t\in[0,\infty): \ \|u^{\lambda}\|_{M_\sigma([0,t]\times\R^n)}\le C_1\lambda^{\frac {s_c[n-3-\sigma(n-6)]}{n-3+4\sigma}}\}.
\end{equation*}
We want to show that $\Omega_1=[0,\infty)$.  We achieve this via a bootstrap argument.  Let
$$
\Omega_2:=\{t\in[0,\infty): \ \|u^{\lambda}\|_{M_\sigma([0,t]\times\R^n)}\le 2C_1\lambda^{\frac {s_c[n-3-\sigma(n-6)]}{n-3+4\sigma}}\}.
$$
In order to run the bootstrap argument successfully, we need to verify four things:
1) $\Omega_1$ is nonempty (as $0\in\Omega_1$),\\
2) $\Omega_1$ is closed (by Fatou's Lemma),\\
3) $\Omega_2\subset\Omega_1$,\\
4) If $T\in \Omega_1$, then there exists $\eps>0$ such that $[T,T+\eps)\subset\Omega_2$; this is a consequence of the local well-posedness theory and
the proof of 3).

We now show 3).  Let $T\in \Omega_2$; we will show that $T\in \Omega_1$.  By \eqref{inter Mora esti} and mass conservation, we have
\begin{align}
\|u^{\lambda}\|_{M_\sigma([0,T]\times\R^n)}
&\lesssim \bigl(\|u_0^\lambda\|_{L_x^2}\|u^\lambda\|_{\hf([0,T]\times\R^n)}\bigr)^{\frac {2\sigma}{n-3+4\sigma}}\|u^\lambda\|_{L_t^\infty\dot H^\sigma_x ([0,T]\times\R^n)}^{\frac{n-3}{n-3+4\sigma}}\nonumber\\
&\lesssim C(\|u_0\|_2)\lambda^{\frac {2\sigma s_c}{n-3+4\sigma}}\|u^{\lambda}\|_{\hf([0,T]\times\R^n)}^{\frac {2\sigma}{n-3+4\sigma}}\|u^\lambda\|_{L_t^\infty\dot H^\sigma_x ([0,T]\times\R^n)}^{\frac{n-3}{n-3+4\sigma}}.\label{Morawetz on 0 T}
\end{align}
To control the second and the third factor, we decompose
$$
u^{\lambda}(t):=P_{\le N}u^{\lambda}(t)+P_{>N}u^{\lambda}(t).
$$
To estimate the low frequencies, we interpolate between the $L_x^2$-norm and $\dot H_x^1$-norm and use the fact that $I$ is the
identity on frequencies $|\xi|\leq N$:
\begin{align}
\|P_{\le N}u^{\lambda}(t)\|_{\dot H^\sigma_x}
&\lesssim \|P_{\le N}u^{\lambda}(t)\|_{L_x^2}^{1-\sigma} \|P_{\le N}u^{\lambda}(t)\|_{\dot H_x^1}^\sigma\notag\\
&\lesssim \lambda^{s_c(1-\sigma)}C(\|u_0\|_2)\|Iu^{\lambda}(t)\|_{\dot H_x^1}^\sigma\\
\|P_{\le N}u^{\lambda}(t)\|_{\dot H^{\frac 12}_x}
&\lesssim \|P_{\le N}u^{\lambda}(t)\|_{L_x^2}^{\frac 12} \|P_{\le N}u^{\lambda}(t)\|_{\dot H_x^1}^{\frac 12}\notag\\
&\lesssim \lambda^{\frac {s_c}2}C(\|u_0\|_2)\|Iu^{\lambda}(t)\|_{\dot H_x^1}^{\frac 12}.\label{low}
\end{align}
To control the high frequencies, we interpolate between $L_x^2$ and $\dot H^s_x$ and use Lemma~\ref{basic property} to get
\begin{align}
\|P_{> N}u^{\lambda}(t)\|_{\dot H_x^\sigma}
&\lesssim \|P_{> N}u^{\lambda}(t)\|_{L_x^2}^{1-\frac \sigma{s}}\|P_{> N}u^{\lambda}(t)\|_{\dot H^s_x}^{\frac \sigma{s}} \notag\\
&\lesssim \lambda^{s_c(1-\frac{\sigma}{s})}\|u_0\|_{L_x^2}^{1-\frac{\sigma}{s}}N^{\frac{\sigma(s-1)}{s}}\|I u^{\lambda}(t)\|_{\dot H_x^1}^{\frac \sigma{s}}\notag\\
&\lesssim \lambda^{s_c(1-\sigma)}C(\|u_0\|_2)\|Iu^{\lambda}(t)\|_{\dot H_x^1}^{\frac \sigma{s}}\\
\|P_{> N}u^{\lambda}(t)\|_{\dot H_x^{\frac 12}}
&\lesssim \|P_{> N}u^{\lambda}(t)\|_{L_x^2}^{1-\frac 1{2s}}\|P_{> N}u^{\lambda}(t)\|_{\dot H^s_x}^{\frac 1{2s}} \notag\\
&\lesssim \lambda^{(1-\frac 1{2s})s_c}N^{\frac {s-1}{2s}}\|u_0\|_{L_x^2}^{1-\frac 1{2s}}\|Iu^{\lambda}(t)\|_{\dot H_x^1}^{\frac 1{2s}}\nonumber\\
&\lesssim \lambda^{\frac {s_c}2}C(\|u_0\|_2)\|Iu^{\lambda}(t)\|_{\dot H_x^1}^{\frac 1{2s}}.\label{high}
\end{align}
Collecting \eqref{Morawetz on 0 T} through \eqref{high}, we obtain
\begin{align}
\|u^{\lambda}\|_{M_\sigma([0,T]\times\R^n)}\lesssim C(&\|u_0\|_2)\lambda^{\frac {s_c[n-3-\sigma(n-6)]}{n-3+4\sigma}}\label{u lambda Morawetz}\notag\\
&\times\sup_{[0,T]}\bigl(\|\nabla I u^{\lambda}(t)\|_2^{\frac 12}+\|\nabla I u^{\lambda}(t)\|_2^{\frac 1{2s}}\bigr)^{\frac {2\sigma}{n-3+4\sigma}}\notag\\
&\times\sup_{[0,T]}\bigl(\|\nabla I u^{\lambda}(t)\|_2^\sigma+\|\nabla I u^{\lambda}(t)\|_2^{\frac{\sigma}{s}}\bigr)^{\frac {n-3}{n-3+4\sigma}}.
\end{align}
Thus, taking $C_1$ sufficiently large depending on $\|u_0\|_{L_x^2}$, we get $T\in \Omega_1$, provided we can establish\footnote{The bootstrap
condition 4) follows from \eqref{bdd kinetic}, Lemma~\ref{basic property}, the local well-posedness theory, \eqref{inter Mora esti}
and the Dominated Convergence Theorem.}
\begin{equation}\label{bdd kinetic}
\sup_{[0,T]}\|\nabla I u^{\lambda}(t)\|_2\le 1.
\end{equation}

We now prove that $T\in \Omega_2$ implies \eqref{bdd kinetic}.  Indeed, let $\eta>0$ be a sufficiently small constant (as in
Proposition~\ref{energy increment scattering}). Divide $[0,T]$ into
$$
L\sim \Bigl(\frac {\lambda^{\frac {s_c[n-3-\sigma(n-6)]}{n-3+4\sigma}}}{\eta}\Bigr)^{\frac{n-3+4\sigma}{\sigma}}\sim\lambda^{\frac{s_c[n-3-\sigma(n-6)]}{\sigma}}
$$
subintervals $I_j=[t_j,t_{j+1}]$, such that for each $j$
$$
\|u^{\lambda}\|_{M_\sigma(I_j\times\R^n)}\leq \eta.
$$
Applying Proposition~\ref{energy increment scattering} on each of the subintervals $I_j$, we get
\begin{equation}\label{kinetic growth}
\sup_{0\le t\le T}\|\nabla I u^{\lambda}(t)\|_2 \le E(Iu_0^\lambda)+C(E(Iu_0^\lambda))\, L\, N^{\min\{1,p\}(s_c-s)+}.
\end{equation}
Then, \eqref{bdd kinetic} follows from \eqref{kinetic growth} as long as
$$
C(E(Iu_0^\lambda))\,L\, N^{\min\{1,p\}(s_c-s)+} \leq E(Iu_0^\lambda)\leq c(\|u_0\|_{H^s_x})  \ll 1.
$$
As $L=O(\lambda^{\frac{s_c[n-3-\sigma(n-6)]}{\sigma}})$, we need to choose $N$ and $\lambda$ such that
$$
N^{1-s}\lambda^{s_c-s}\|u_0\|_{H^s_x}\ll1 \quad \text{and} \quad \lambda^{\frac{s_c[n-3-\sigma(n-6)]}{\sigma}}N^{\min\{1,p\}(s_c-s)+}\leq c(\|u_0\|_{H^s_x})\ll1.
$$
Plugging the first relation into the second one, we see that we need to choose $N$ depending on $\|u_0\|_{H^s_x}$ such that
$$
N^{\frac {s_c(1-s)[n-3-\sigma(n-6)]}{\sigma(s-s_c)}+\min\{1,p\}(s_c-s)+}\ll1.
$$
This is possible whenever $s$ is such that
\begin{align*}
s_c(1-s)[n-3-\sigma(n-6)]<\min\{1,p\}\sigma(s_c-s)^2,
\end{align*}
i.e., $s>s_+(n,p,\sigma)$ where $s_+(n,p,\sigma)$ is the larger of the two roots to the quadratic equation
\begin{align*}
s_c(1-s)[n-3-\sigma(n-6)]=\min\{1,p\}\sigma(s_c-s)^2.
\end{align*}

Thus, the bootstrap is complete and \eqref{rescaled Mora} follows.  Hence, \eqref{bdd kinetic} holds for all $T\in \R$, which by \eqref{i3}
and the conservation of mass implies
\begin{align*}
\|u(T)\|_{H_x^s}
&\lesssim \|u_0\|_{L_x^2} + \|u(T)\|_{\dot H^s_x}\\
&\lesssim \|u_0\|_{L_x^2} + \lambda^{s-s_c}\|u^\lambda(\lambda^2T)\|_{\dot H^s_x}\\
&\lesssim \|u_0\|_{L_x^2} + \lambda^{s-s_c}\|Iu^\lambda(\lambda^2T)\|_{H^1_x}\\
&\lesssim \|u_0\|_{L_x^2} + \lambda^{s-s_c}\bigl(\|u^\lambda(\lambda^2T)\|_{L^2_x}+\|\nabla Iu^\lambda(\lambda^2T)\|_{L_x^2}\bigr)\\
&\lesssim \|u_0\|_{L_x^2} + \lambda^{s-s_c}(\lambda^{s_c}\|u_0\|_{L_x^2}+1)\\
&\leq C(\|u_0\|_{H^s_x}),
\end{align*}
for all $T\in \R$.  Therefore,
\begin{align}\label{global H^s}
\|u\|_{L_t^\infty H^s_x}\leq C(\|u_0\|_{H^s_x}).
\end{align}

\subsection{Scattering}
We first show that the global Morawetz estimate \eqref{global Morawetz} can be upgraded to the global Strichartz bound
\begin{equation}\label{w bound}
\|u\|_W:=\sup_{(q,r) admissible}\|\langle\nabla\rangle^s u\|_{L_t^qL_x^r(\R\times\R^n)}\leq C(\|u_0\|_{H^s_x}).
\end{equation}
The second step is to use this estimate to prove asymptotic completeness. The construction of the wave operators is standard and we omit it.

Let $u$ be a global solution to \eqref{equation} with initial data in $H^s(\R^n)$ for $s>s_0(n,p)$.  By \eqref{global Morawetz} we have
$$
\|u\|_{M_\sigma(\R\times\R^n)}\le C(\|u_0\|_{H^s_x}).
$$
Let $\delta>0$ be a small constant to be chosen momentarily and split $\R$ into $L=O(\|u_0\|_{H_x^s})$ subintervals $I_j=[t_j,t_{j+1}]$ such that
$$
\|u\|_{M_\sigma(I_j\times\R^n)}\le \delta.
$$
By Strichartz,
\begin{equation}\label{wuij}
\|u\|_{W(I_j)}
\lesssim\|\langle\nabla\rangle^su(t_j)\|_2+\|\langle\nabla\rangle^s(|u|^pu)\|_{L_t^2L_x^{\frac {2n}{n+2}}(I_j\times\R^n)}.
\end{equation}
Using Lemma~\ref{use Morawetz} and \eqref{global H^s}, we control the nonlinearity as follows:
\begin{align}
\|\langle\nabla\rangle^s(|u|^pu)\|_{L_t^2L_x^{\frac {2n}{n+2}}(I_j\times\R^n)}
&\lesssim \|u\|_{W(I_j)}\|u\|_{M_\sigma(I_j\times \R^n)}^{\frac{\eps(n-3+4\sigma)}{2\sigma(2+\eps)}}\|u\|^{\alpha(\eps)+\beta(\eps)}_{L_t^{\infty} H_x^s(I_j\times \R^n)}\notag\\
&\lesssim \|u\|_{W(I_j)}\delta^{\frac{\eps(n-3+4\sigma)}{2\sigma(2+\eps)}} C(\|u_0\|_{H^s_x}).\label{nonlinearity}
\end{align}
Taking $\delta$ sufficiently small depending only on $\|u_0\|_{H^s_x}$, \eqref{wuij} and \eqref{nonlinearity} yield
$$
\|u\|_{W(I_j)}\lesssim \|\langle\nabla\rangle^su(t_j)\|_2.
$$
Adding these bounds over all subintervals $I_j$, we obtain \eqref{w bound}.

We now use \eqref{w bound} to show asymptotic completeness, i.e., there exist unique $u_{\pm}\in H^s_x$ such that
$$
\lim_{t\to \pm\infty}\|u(t)-e^{it\Delta}u_{\pm}\|_{H^s_x}=0.
$$
By time reversal symmetry, it suffices to argue in the positive time direction.  For $t>0$ define $v(t) = e^{-it\Delta}u(t)$. We will show
that $v(t)$ converges in $H^s_x$ as $t\rightarrow \infty$, and define $u_+$ to be that limit.

Indeed, from Duhamel's formula \eqref{duhamel} we have
\begin{align}\label{v}
v(t) = u_0 - i\int_{0}^{t} e^{-is\Delta}(|u|^pu)(s)ds.
\end{align}
Therefore, for $0<\tau<t$,
$$
v(t)-v(\tau)=-i\int_{\tau}^{t}e^{-is\Delta}(|u|^pu)(s)ds.
$$
By Strichartz and Lemma~\ref{use Morawetz}, we estimate
\begin{align*}
\|v(t)-v(\tau)\|_{H^s_x}
&=\|e^{it\Delta}[v(t)-v(\tau)]\|_{H^s_x} \\
&\lesssim \|\langle\nabla\rangle^s(|u|^pu)\|_{L_t^2L_x^{\frac {2n}{n+2}}([\tau,t]\times\R^n)}\\
&\lesssim \|u\|_{W([\tau,t])}\|u\|_{M_\sigma([\tau,t]\times\R^n)}^{\frac{\eps(n-3+4\sigma)}{2\sigma(2+\eps)}}\|u\|^{\alpha(\eps)+\beta(\eps)}_{L_t^{\infty} H_x^s([\tau,t]\times\R^n)}.
\end{align*}
Using \eqref{global Morawetz}, \eqref{global H^s}, and \eqref{w bound}, we obtain
\begin{center}
$\|v(t)-v(\tau)\|_{H^s_x}\rightarrow 0 \quad$ as $\tau,t \rightarrow \infty$.
\end{center}
In particular, this implies that $u_{+}$ is well defined. Also, inspecting \eqref{v} one easily sees that
\begin{align}
u_{+}=u_0- i\int_{0}^{\infty}e^{-is\Delta}(|u|^pu)(s)ds
\end{align}
and thus
\begin{align}\label{u+}
e^{it\Delta}u_{+}=e^{it\Delta}u_0- i\int_{0}^{\infty}e^{i(t-s)\Delta}(|u|^pu)(s)ds.
\end{align}
By the same arguments as above, \eqref{u+} and Duhamel's formula \eqref{duhamel} imply that
$\|u(t)-e^{it\Delta}u_{+}\|_{H^s_x}\rightarrow 0$ as $t\rightarrow\infty$.

\subsection{Proof of Proposition~\ref{energy increment scattering}}
The first step is to upgrade Lemma~\ref{nablaif}.  We are interested in controlling commutators in spacetime norms, not merely pointwise in time.

For any spacetime slab $\ir$, we define
$$
Z_I:=\sup_{(q,r) admissible}\|\nabla Iu\|_{L_t^qL_x^r(\ir)}.
$$

\begin{lemma}\label{gif scat}
Let $I$ be a compact time interval,
$\frac{1+\min\{1,p\}s_c}{1+\min\{1,p\}}<s<1$, $0<\sigma\leq s$, and
$\tfrac {4(n-3-4\sigma)}{n(n-3+2\sigma)+4\sigma}<p<\tfrac{4}{n-2s}$.
Assume that
$$
\|u\|_{M_\sigma(\ir)}\leq \eta,
$$
for a small constant $\eta>0$.  Then,
\begin{align}
&\|\nabla I F(u)-(I\nabla u)F'(u)\|_{\eddp} \lesssim N^{\min\{1,p\}(s_c-s)+}Z_I\bigl(\eta^{\theta}Z_I^{1-\theta} + Z_I\bigr)^p\label{gif scat1}\\
&\|\nabla I F(u)\|_{\eddp}\lesssim N^{\min\{1,p\}(s_c-s)+}Z_I^{p+1} +
\eta^{p\theta}Z_I^{1+p(1-\theta)}.\label{gif scat2}
\end{align}
Here, $\theta$ is as defined in Lemma~\ref{interpolation1}.
\end{lemma}

\begin{proof}
Throughout the proof all spacetime norms will be on $\ir$.

As by assumption $\frac{1+\min\{1,p\}s_c}{1+\min\{1,p\}}<s$, there exists $\delta_0>0$ such that for any $0<\delta<\delta_0$ we have
$\frac{1+\min\{1,p\}(s_c+\delta)}{1+\min\{1,p\}}<s$.  Let $\nu:=-1+s+\min\{1,p\}(s-s_c-\delta)$; it is easy to check that we have $0<\nu<s$.
Applying Lemma~\ref{nablaif} with this value of $\nu$ and using H\"older in time, we obtain
\begin{align*}
\|\nabla &I F(u)-(I\nabla u)F'(u)\|_{\eddp}\\
&\lesssim N^{-\min\{1,p\}(s-s_c-\delta)}\|\nabla Iu\|_{\ess}\|\lnr^{\min\{1,p\}(s-s_c-\delta)}F'(u)\|_{\er}\\
&\lesssim N^{\min\{1,p\}(s_c-s+\delta)}Z_I\|\lnr^{\min\{1,p\}(s-s_c-\delta)}F'(u)\|_{\er}.
\end{align*}
Here, $\eps$ is as in Lemma~\ref{interpolation1}.  The estimate \eqref{gif scat1} follows from the above estimate, provided
\begin{align}\label{commut show1}
\|\lnr^{\min\{1,p\}(s-s_c-\delta)}F'(u)\|_{\er}\lesssim \bigl(\eta^{\theta}Z_I^{1-\theta} + Z_I\bigr)^p.
\end{align}

If $\min\{1,p\}=1$, we bound
\begin{align*}
\|\lnr^{\min\{1,p\}(s-s_c-\delta)}&F'(u)\|_{\er}\\
&\lesssim\|\lnr^{\min\{1,p\}(s-s_c)}F'(u)\|_{\er}\\
&\lesssim\|F'(u)\|_{\er}+\||\nabla|^{\min\{1,p\}(s-s_c)}F'(u)\|_{\er}.
\end{align*}
Clearly,
$$
\|F'(u)\|_{\er}\lesssim\|u\|_{\ep}^p.
$$
Decomposing $u:=u_{\leq N}+u_{>N}$, by hypothesis, \eqref{i2},
Lemma~\ref{interpolation1}, Sobolev embedding, and the fact that the
operator $I$ is the identity on frequencies $|\xi|\leq N$, we
estimate
\begin{align}
\|u_{\le N}\|_{\ep} &\lesssim  \|u_{\le
N}\|_{M_\sigma}^{\theta}\|\nabla u_{\le
N}\|_{\frac{n-3+4\sigma}{\sigma},\frac{2n(n-3+4\sigma)}{n(n-3+4\sigma)-4\sigma}}^{1-\theta}
  \lesssim \eta^{\theta}Z_I^{1-\theta}\label{ulo commut}\\
\|u_{>N}\|_{\ep} &\lesssim \||\nabla|^{s_c} u_{>N}\|_{\epp}
  \lesssim N^{s_c-1} Z_I.\label{uhigh}
\end{align}
Hence, for $N$ sufficiently large,
\begin{equation}\label{estimate of u}
\|u\|_{\ep}\lesssim \eta^{\theta}Z_I^{1-\theta} +Z_I
\end{equation}
and thus,
\begin{align}\label{commut show2}
\|F'(u)\|_{\er}\lesssim \bigl(\eta^{\theta}Z_I^{1-\theta} +
Z_I\bigr)^p.
\end{align}
Using Lemma~\ref{F Lip}, H\"older in time, and \eqref{estimate of u}, we estimate
\begin{align}
\||\nabla|^{\min\{1,p\}(s-s_c)}&F'(u)\|_{\er}\notag\\
&=\||\nabla|^{s-s_c}F'(u)\|_{\er}\notag\\
&\lesssim \||\nabla|^{s-s_c}u\|_{\frac{2p(2+\eps)}{\eps}, \frac{np(2+\eps)}{4+\eps}}\|u\|^{p-1}_{\frac{2p(2+\eps)}{\eps}, \frac{np(2+\eps)}{4+\eps}}\notag\\
&\lesssim \||\nabla|^{s-s_c}u\|_{\frac{2p(2+\eps)}{\eps},\frac{np(2+\eps)}{4+\eps}}\bigl(\eta^{\theta}Z_I^{1-\theta} +Z_I\bigr)^{p-1}.\label{p large}
\end{align}

If instead $\min\{1,p\}=p$, we bound
\begin{align*}
\|&\lnr^{\min\{1,p\}(s-s_c-\delta)}F'(u)\|_{\er}\\
&\qquad \lesssim\|F'(u)\|_{\er}+\||\nabla|^{\min\{1,p\}(s-s_c-\delta)}F'(u)\|_{\er}.
\end{align*}
Using Lemma~\ref{fdfp} (with $\alpha:=p$, $\sigma:=p(s-s_c)$, and $\rho:=s-s_c$), H\"older in time, and \eqref{estimate of u}, we get
\begin{align}
\||\nabla|^{\min\{1,p\}(s-s_c-\delta)}&F'(u)\|_{\er}\notag\\
&=\||\nabla|^{p(s-s_c-\delta)}F'(u)\|_{\er}\notag\\
&\lesssim \|u\|_{\ep}^{\frac{p\delta}{s-s_c}}\||\nabla|^{s-s_c}u\|_{\ep}^{\frac{p(s-s_c-\delta)}{s-s_c}}\notag\\
&\lesssim \bigl(\eta^{\theta}Z_I^{1-\theta} +Z_I\bigr)^{\frac{p\delta}{s-s_c}}\||\nabla|^{s-s_c}u\|_{\ep}^{\frac{p(s-s_c-\delta)}{s-s_c}}.\label{psmall}
\end{align}

To estimate $\||\nabla|^{s-s_c}u\|_{\ep}$, we decompose $u:=u_{\le 1}+u_{1<\cdot\le N}+u_{>N}$. The very low frequencies we control by
Lemma~\ref{interpolation1}, Bernstein, and the fact that $I$ is the identity on frequencies $|\xi|\leq N$:
\begin{align*}
\||\nabla|^{s-s_c} u_{\le 1}\|_{\ep}
&\lesssim \|u_{\le 1}\|_{\ep}\\
&\lesssim \|u_{\le1}\|_{M_\sigma}^{\theta}\|\nabla u_{\le1}\|_{\frac{n-3+4\sigma}{\sigma},\frac{2n(n-3+4\sigma)}{n(n-3+4\sigma)-4\sigma}}^{1-\theta}\\
&\lesssim \eta^{\theta} Z_I^{1-\theta}.
\end{align*}
For the medium frequencies, we use Sobolev embedding, Bernstein, and the definition of $I$ to obtain
\begin{align*}
\||\nabla|^{s-s_c} u_{1\le \cdot\le N}\|_{\ep}
&\lesssim \||\nabla|^s u_{1<\cdot\le N}\|_{\epp}\\
&\lesssim \|\nabla Iu_{1<\cdot\le N}\|_{\epp}\\
&\lesssim Z_I.
\end{align*}
To estimate the high frequencies, we use Sobolev embedding and Lemma~\ref{basic property}:
\begin{align*}
\||\nabla| ^{s-s_c} u_{>N}\|_{\ep}
&\lesssim \||\nabla|^s u_{>N}\|_{\epp}\lesssim N^{s-1} Z_I.
\end{align*}
As $s-1<0$, for $N$ sufficiently large we get
\begin{equation}\label{deriv of u}
\||\nabla|^{s-s_c} u\|_{\ep}\lesssim \eta^{\theta}Z_I^{1-\theta} + Z_I.
\end{equation}

Collecting \eqref{commut show2}, \eqref{p large}, \eqref{psmall}, and \eqref{deriv of u}, we derive \eqref{commut show1} and hence
the claim \eqref{gif scat1}.

We turn now toward \eqref{gif scat2}.  By H\"older, \eqref{ulo
commut}, and \eqref{uhigh},
\begin{align*}
\|(I\nabla u)F'(u)\|_{\eddp}
&\lesssim \|\nabla Iu\|_{\ess}\|F'(u)\|_{\er}\\
&\lesssim Z_I\bigl(\eta^{\theta}Z_I^{1-\theta}+ N^{s_c-1}Z_I\bigr)^p\\
&\lesssim \eta^{p\theta} Z_I^{1+p(1-\theta)} +N^{p(s_c-1)}Z_I^{1+p}\\
&\lesssim \eta^{p\theta} Z_I^{1+p(1-\theta)}+N^{\min\{1,p\}(s_c-s)+}Z_I^{1+p}.
\end{align*}
This, together with \eqref{gif scat1} and the triangle inequality imply \eqref{gif scat2}.
\end{proof}

\vspace{0.4cm}

We are now ready to prove Proposition~\ref{energy increment
scattering}. The proof is carried out in two steps.  First, we use
Strichartz estimates and the Morawetz control to estimate $Z_I.$
Indeed, by Strichartz and Lemma~\ref{gif scat} (here, $I:=[t_0,T]$),
we have
\begin{align*}
Z_I
&\lesssim \|\nabla Iu(t_0)\|_2+\|\nabla IF(u)\|_{\eddp}\\
&\lesssim \|\nabla Iu(t_0)\|_2+N^{\min\{1,p\}(s_c-s)+}Z_I^{1+p}
+\eta^{p\theta} Z_I^{1+p(1-\theta)}.
\end{align*}
As by assumption $\|\nabla Iu(t_0)\|_2\leq 1$, a standard continuity argument yields
\begin{equation}\label{zit control}
Z_I\lesssim \|\nabla Iu(t_0)\|_2,
\end{equation}
provided $\eta>0$ is chosen sufficiently small and $N$ is chosen sufficiently large.

Next, we express the energy increment in terms of $Z_I$. In what follows, all spacetime norms are taken on the slab $\ir:=[t_0,T]\times\R^n$.

By the Fundamental Theorem of Calculus,
\begin{align*}
E(Iu(T))-E(Iu(t_0))&=\int_{t_0}^T \frac{\partial}{\partial s}E(Iu(s)) \,ds\\
                    &=\Re \int_{t_0}^T \int_{\R^n} \overline{Iu_t}(-\Delta Iu+F(Iu)) \,dx\,ds.
\end{align*}
As $Iu_t=i\Delta Iu-iIF(u)$, we have
$$
\Re \int_{t_0}^T\int_{\R^n}\overline {Iu_t} (-\Delta Iu+IF(u)) \,dx\,ds=0,
$$
and so, after an integration by parts, we can write the energy increment as
\begin{align}
E(Iu(T))-E(Iu(t_0))
&=-\Im \int_{t_0}^T\int_{\R^n}\overline{\nabla Iu}\nabla(F(Iu)-IF(u))\,dx\,dt\nonumber\\
&\quad-\Im \int_{t_0}^T\int_{\R^n}\overline{IF(u)}(F(Iu)-IF(u))\,dx\,dt\nonumber\\
&=-\Im \int_{t_0}^T\int_{\R^n}\overline{\nabla Iu}\cdot \nabla Iu(F'(Iu)-F'(u))\,dx\,dt\label{term1}\\
&\quad-\Im\int_{t_0}^T\int_{\R^n}\overline{\nabla Iu}\cdot((\nabla Iu)F'(u)-I(F'(u)\nabla u))\,dx\,dt\label{term2}\\
&\quad-\Im \int_{t_0}^T\int_{\R^n}\overline{IF(u)}(F(Iu)-IF(u))\,dx\,dt\label{term3}.
\end{align}

Consider \eqref{term1}. By H\"older and \eqref{holder continuity}, we get
\begin{align*}
|\eqref{term1}|
&\lesssim \|\nabla Iu\|_{2,\frac{2n}{n-2}}\|\nabla Iu\|_{\ess} \|F'(Iu)-F'(u)\|_{\er} \\
&\lesssim Z_I^2\||Iu-u|^{\min\{1,p\}}(|u|+|Iu|)^{p-\min\{1,p\}}\|_{\er}\\
&\lesssim Z_I^2\|P_{> N}u\|_{\ep}^{\min\{1,p\}}\|u\|_{\ep}^{p-\min\{1,p\}}.
\end{align*}
By \eqref{uhigh} and \eqref{estimate of u}, we obtain
\begin{align}
|\eqref{term1}|
&\lesssim N^{\min\{1,p\}(s_c-1)} Z_I^2 Z_I^{\min\{1,p\}} \bigl(\eta^\theta Z_I^{1-\theta}+Z_I\bigr)^{p-\min\{1,p\}}\notag\\
&\lesssim N^{\min\{1,p\}(s_c-s)+}Z_I^2\bigl(\eta^\theta
Z_I^{1-\theta}+Z_I\bigr)^p.\label{control term1}
\end{align}

To estimate \eqref{term2} we use H\"older and \eqref{gif scat1}:
\begin{align}
|\eqref{term2}|&\lesssim \|\nabla Iu\|_{2,\frac {2n}{n-2}}\|(\nabla Iu)F'(u)-I(F'(u)\nabla u)\|_{2,\frac {2n}{n+2}}\notag\\
&\lesssim N^{\min\{1,p\}(s_c-s)+}Z_I^2\bigl(\eta^\theta
Z_I^{1-\theta}+Z_I\bigr)^p.\label{control term2}
\end{align}

We now consider \eqref{term3}. By H\"older and Sobolev embedding, we estimate
\begin{align*}
|\eqref{term3}|
&\lesssim \||\nabla|^{-1}I F(u)\|_{2,\frac {2n}{n-2}}\|\nabla[F(Iu)-IF(u)]\|_{2,\frac {2n}{n+2}}\\
&\lesssim \|\nabla IF(u)\|_{2,\frac{2n}{n+2}}\|\nabla[F(Iu)-IF(u)]\|_{2,\frac {2n}{n+2}}.
\end{align*}
By \eqref{gif scat2},
\begin{align*}
\|\nabla IF(u)\|_{2,\frac{2n}{n+2}}
&\lesssim N^{\min\{1,p\}(s_c-s)+}Z_I^{1+p} +\eta^{p\theta}Z_I^{1+p(1-\theta)}\\
&\lesssim Z_I\bigl(Z_I+\eta^\theta Z_I^{1-\theta}\bigr)^p.
\end{align*}
Meanwhile, by the triangle inequality,
\begin{align*}
\|\nabla[F(Iu)-IF(u)]\|_{2,\frac {2n}{n+2}}
&\lesssim \|(\nabla Iu)[F'(Iu)-F'(u)]\|_{2,\frac {2n}{n+2}}\\
&\quad +\|(\nabla Iu)F'(u)-I(\nabla u\cdot F'(u))\|_{2,\frac {2n}{n+2}}.
\end{align*}
Using the same argument as that used to estimate \eqref{term1}, we get
$$
\|(\nabla Iu)[F'(Iu)-F'(u)]\|_{2,\frac {2n}{n+2}}
\lesssim N^{\min\{1,p\}(s_c-s)+}Z_I\bigl(\eta^\theta Z_I^{1-\theta}+Z_I\bigr)^p,
$$
while by \eqref{gif scat1}, we have
$$
\|(\nabla Iu)F'(u)-I(\nabla u\cdot F'(u))\|_{2,\frac {2n}{n+2}}\lesssim N^{\min\{1,p\}(s_c-s)+}Z_I\bigl(\eta^\theta Z_I^{1-\theta}+Z_I\bigr)^p.
$$
Thus,
\begin{align}
|\eqref{term3}|\lesssim N^{\min\{1,p\}(s_c-s)+}Z_I^2\bigl(\eta^\theta Z_I^{1-\theta}+Z_I\bigr)^{2p}.\label{control term3}
\end{align}

Putting together \eqref{control term1}, \eqref{control term2}, and \eqref{control term3}, and recalling \eqref{zit control},
we obtain the desired bound on the energy increment, i.e.,
\begin{align*}
|E(Iu(T))-E(Iu(t_0))|
&\lesssim N^{\min\{1,p\}(s_c-s)+}Z_I^2\bigl[(\eta^\theta Z_I^{1-\theta}+Z_I\bigr)^p+(\eta^\theta Z_I ^{1-\theta}+Z_I\bigr)^{2p}\bigr]\\
&\leq C(\|\nabla I u(t_0)\|_2)N^{\min\{1,p\}(s_c-s)+}.
\end{align*}

%%%%%%%%%%%%%%%%%%%%%%%%%%%%%%%%%%%%%%%%%%%%%%%%%%%%%%%%%%%%%%%%%%%%%%%%%%%%%%%%%%%%%%%%%%%
%
%
%                                   Section
%
%
%%%%%%%%%%%%%%%%%%%%%%%%%%%%%%%%%%%%%%%%%%%%%%%%%%%%%%%%%%%%%%%%%%%%%%%%%%%%%%%%%%%%%%%%%%%

\end{document}